\input amstex
\documentstyle{amsppt}\nologo\footline={}

\def\Homeo{\mathop{\text{\rm Homeo}}}
\def\Isom{\mathop{\text{\rm Isom}}}
\def\Cl{\mathop{\text{\rm Cl}}}
\def\T{\mathop{\text{\rm T}}}
\def\N{\mathop{\text{\rm N}}}

\input epsf
\input lpic-mod

\hsize450pt\vsize591pt\topmatter\title Poincar\'e's polyhedron theorem for cocompact groups in dimension $4$\endtitle\thanks
Second author partially supported by FAPESP (grant 2012/07587-4) and CNPq (grant 304671/2012-6).\endthanks\author Sasha Anan$'$in,
Carlos H.~Grossi, J\'ulio C.~C.~da Silva\endauthor\address Departamento de Matem\'atica, IMECC, Universidade Estadual de
Campinas,\newline13083-970--Campinas--SP, Brasil\endaddress\address Departamento de Matem\'atica, ICMC, Universidade de S\~ao
Paulo, Caixa Postal 668, \newline13560-970--S\~ao Carlos--SP, Brasil\endaddress\address Departamento de Matem\'atica, IMECC,
Universidade Estadual de Campinas,\newline13083-970--Campinas--SP, Brasil\endaddress\abstract
We prove a version of Poincar\'e's polyhedron theorem whose requirements are as local as possible. New techniques such as the use
of discrete groupoids of isometries are introduced. The theorem may have a wide range of applications and can be generalized to the
case of higher dimension and other geometric structures. It is planned as a first step in a program of constructing compact
$\Bbb C$-surfaces of general type satisfying $c_1^2=3c_2$.\endabstract\endtopmatter\document

{\hfill\it To Misha Verbitsky on no particular occasion}

\bigskip

\bigskip

\centerline{\bf1.~Introduction}

\medskip

{\bf1.1.}~One of the few powerful tools for constructing compact manifolds (or orbifolds) carrying a prescribed geometry is closely
related to Poincar\'e's polyhedron theorem (PPT). A compact manifold $X$ in question is the quotient of a model space $M$ by a
discrete cocompact group $G$ of automorphisms. One can associate to $G$ a polyhedron $P\subset M$ into which $X=M/G$ can be `cut
and unfolded.' For instance, at the presence of metric, there is the compact Dirichlet polyhedron
$P_c:=\big\{x\in M\mid\forall\ g\in G\ d(x,c)\le d(x,gc)\big\}$ centred at a generic point $c\in M$; it has finitely many
codimension $1$ faces (contained in hypersurfaces equidistant from two points) and face-pairing isometries between such faces. In
general, the settings of PPT include a polyhedron whose face-pairing automorphisms generate $G$ and, {\it a posteriori,} the
discreteness of $G$ turns out to be a consequence of the tessellation of $M$ by the $G$-copies of the polyhedron. In this way, PPT
transforms constructing compact manifolds into verifying tessellation.

Generally speaking, the tessellation of $M$ resides in the mutual position of certain $G$-copies of codimension $1$ faces of $P$.
In the case of classic geometries [AGr1], the hypersurfaces equidistant from two points are real algebraic and rather simple. They
constitute a reasonable choice for codimension $1$ faces of $P$. Nevertheless, inferring the tessellation of $M$ from an analysis
of the mutual position of some $G$-copies of codimension $1$ faces is frequently impossible to perform. In this paper, we follow
the strategy taken from [AGr2] and develop a version of PPT whose requirements providing the tessellation are as local as possible,
hence, easily verifiable in practice. This allows to avoid the mentioned difficulty.

\medskip

{\bf1.2.}~In few vague words, Theorems 4.1 and 4.3 claim that the following requirements provide a tessellation of $M$. A small
neighbourhood of a generic point of every codimension $1$ or $2$ face of $P$ is required to be tessellated by `expected' copies of
$P$. Another condition asks the polyhedron to remain locally connected after removing its codimension $3$ faces. The last
requirement (empty in the $3$-dimensional case) is that the stabilizer in $G$ of every codimension $3$ face induces on the face a
finite group of isometries, where $G$ stands for the group generated by the face-pairing isometries.

Local tessellation around a generic point in a face of codimension $1$ or $2$ is translated in Lemmas~3.10 and 3.12 into
infinitesimal conditions called `interior into exterior' and `total angle' conditions; the latter takes the form `sum of angles
$=2\pi$.'

Speaking of tessellation around generic points in codimension $3$ faces, one can formulate in a similar fashion a `solid angle'
condition. The readily verifiable connectedness near codimension $3$ is a topological form of such a `solid angle' condition.
Example 3.16 shows that PPT fails without this connectedness requirement.\footnote{We are grateful to the referee who stimulated us
to state the conjecture in [AGr2, Subsection 4.3]; unfortunately, our conjecture lacks the condition of connectedness near
codimension $3$.}

By definition, the discreteness of $G$ means that $gK\cap K\ne\varnothing$ only for finitely many $g\in G$ if $K\subset M$ is a
compact set. Taking for $K$ a codimension $3$ face we conclude the necessity of the requirement concerning the stabilizer. In high
dimension, this sort of condition does not seem verifiable in general. Already in dimension 4, it falls into the verifiable
category only because those codimension $3$ faces that may cause a problem have to be round circles (allegedly, a rare case). So,
it is hopeless to dream of a general version of PPT (with easily verifiable hypothesis) for dimension $\ge5$.

In literature, there are many versions of PPT requiring that the polyhedron modulo face-pairing is equipped with a naturally
induced complete metric; see, for instance, [Mas, Condition H.~10, page~75]. In~our settings, completeness would be automatic as
polyhedra are compact, but the induced pseudo-metric is not always a metric (the quotient can even be non-Hausdorff). Verifying the
induced pseudo-metric to be a metric is not in general a task that can be performed in practice.

\medskip

{\bf1.3.}~In section 3, we introduce and study a concept of smooth compact polyhedra. On one hand, this concept allows us to
proceed by induction on dimension. On the other hand, it suits all possible practical needs.

A crucial feature of the theorems is the extensive use of groupoids caused by the necessity of dealing with nonconnected spaces and
polyhedra in the proofs. The lack of sufficient symmetries makes it impossible to apply induction on dimension and to keep
simultaneously the symmetries at the level of a group of isometries of the ambient space. Curiously, adapting [AGr2, Section 2] to
the case of groupoids required only a few changes.

It does not seem to be difficult to generalize the results of this paper to higher dimension just by following the presented proofs
(though this may not be as easy as simply writing this phrase). Our interest lies in dimension $4$ : this paper is intended as a
first step in the program [Ana] of constructing compact $\Bbb C$-surfaces of general type satisfying the equality in
Bogomolov-Miyaoka-Yau's inequality $c_1^2\le\nomathbreak3c_2$ between Chern numbers. Such surfaces are compact complex hyperbolic
manifolds [Yau] and attract considerable interest [Rei]. They are hard to construct. The few known examples include those in [Mos]
and the fake projective planes [Mum], [CaS].

\medskip

{\bf1.4.}~The scope of possible applications of the theorems lies far beyond the originally planned one. As far as compact
manifolds are concerned, the ideas involved in the proofs seem to be adaptable to constructing manifolds carrying other geometric
structures. For instance, it seems possible to modify them and incorporate the case of compact affine manifolds.

Theorems 4.1 and 4.3 are not free from the defect indicated at the end of [AGr2, Subsection 1.3]. It is still necessary to verify
the simplicity of the polyhedron, i.e., to study the mutual position of its codimension $1$ faces. Nevertheless, there is some
progress in the direction of getting rid of the global requirement of simplicity. Groupoids allow us to work effectively with
nonconnected polyhedra. So,~we~can partially escape the verification that the faces intersect properly by {\it a priori\/} cutting
the {\it a posteriori\/} problematic polyhedron into small pieces and placing them distantly by means of arbitrary isometries.

\bigskip

\centerline{\bf2.~Discrete groupoids of isometries}

\medskip

In this section, we reformulate the corresponding [AGr2, Section 2] in a bit more general form as we need to drop the assumption
that $M$ or $P$ is connected.

Let $M$ be a locally path-connected metric space with finitely many connected components, denoted by $M_i$, $i\in I$ (in
applications, each $M_i$ will be simply-connected). For convenience, we assume that $d(x_1,x_2)=\infty$ for points $x_1,x_2\in M$
from different components. Denote by $B(x,\varepsilon)$ the open ball of radius $\varepsilon>0$ centred at $x$ and let
$N(X,\varepsilon):=\bigcup\limits_{x\in X}B(x,\varepsilon)$ for $X\subset M$. In this section, we regard a {\it polyhedron\/} in
$M$ as being a closed, locally path-connected subspace $P\subset M$ such that

\smallskip

$\bullet$ $P$ has finitely many connected components;

\smallskip

$\bullet$ $P$ is the closure of its nonempty interior: $\overset{\,_\circ}\to P\ne\varnothing$ and $P=\Cl\overset{\,_\circ}\to P$;

\smallskip

$\bullet$ the nonempty boundary of $P$ is decomposed into the union of nonempty path-connected subsets $s\in S$ called
{\it faces\/}: $\partial P:=P\setminus\overset{\,_\circ}\to P=\bigcup\limits_{s\in S}s$.

\smallskip

The {\it isometry groupoid\/} $\Isom M$ of $M$ is the category whose objects are $M_i$, $i\in I$, and morphisms are all possible
isometries between the $M_i$'s. (By definition, a {\it groupoid\/} is a category whose morphisms are isomorphisms.) When writing
$g\in\Isom M$, we mean that $g$ is a morphism whose domain and codomain are therefore prescribed, i.e., $g:M_i\to M_j$ is an
isometry for suitable $M_i,M_j$.

A {\it face-pairing\/} of a polyhedron $P$ is an involution $\overline{\phantom{m}}:S\to S$ and a family of isometries
$I_s\in\Isom M$ satisfying $I_ss=\overline s$ and $I_{\overline s}=I_s^{-1}$ for every face $s\in S$. Since the faces
$s,\overline s$ are connected, they determine the corresponding components $M_i,M_j$, hence, we have $I_s:M_i\to M_j$. We denote
$P_i:=P\cap M_i$. Note that $P_i$ is not necessarily connected.

We need a brief elementary introduction to groupoids. Groupoids in this section will always have their objects indexed by $I$ and
all functors will be identical at the level of objects. Frequently, the objects in a groupoid will be topological spaces. In this
case, we think of the disjoint union of these objects as of an $I$-{\it graded\/} topological space. Any $I$-graded topological
space $T$ gives raise to a groupoid $\Homeo T$ of its homeomorphisms. An action of a groupoid $G$ on an $I$-graded space $T$ (by
homeomorphisms) is simply a functor $G\to\Homeo T$. In this case, we say that $T$ is a $G$-{\it space.} In order to define an
action of $G$ on~$T$ (by homeomorphisms), one only needs to associate to every $G\ni g:i\to j$ a homeomorphism $T_i\to T_j$ (we
will denote it by the same letter $g:T_i\to T_j$) so that the composition in $G$ provides the composition of the associated
homeomorphisms. (By convention, when we write $gt$, we silently require that $g\in G$ is applicable to $t\in T$.) For example,
$\Isom M$ acts on $M$ (by isometries, hence, by homeomorphisms) due to the inclusion $\Isom M\subset\Homeo M$. Given $I$-graded
$G$-spaces $T,T'$, an $I$-graded continuous map $f:T\to T'$ is said to be a $G$-{\it map\/} if $f(gt)=gf(t)$ for all $g\in G$ and
$t\in T$.

\smallskip

Let $P$ be a polyhedron with a given face-pairing and let $G\subset\Isom M$ denote the subcategory generated by the face-pairing
isometries. Assuming the discrete topology on $G$, we equip $P$, $G\times P$, and $G*P:=\bigsqcup_iG\times_iP\subset G\times P$
with their natural topologies, where $G\times_iP:=\big\{(g,x)\mid x\in P_j,\ G\ni g:M_j\to M_i\big\}$. The rule $(g,x)\mapsto gx$
defines an $I$-graded continuous map $\psi:G*P\to M$, i.e., $\psi(G\times_iP)\subset M_i$.

We introduce a relation in $G*P$ by putting $(g,x)\sim(h,y)$ exactly when $x\in s$ for some $s\in S$, $I_sx=y$, and $h^{-1}g=I_s$.
Taking the closure of this symmetric relation with respect to transitivity (and reflexivity), we obtain an equivalence relation
also denoted by $\sim$. Let $J:=G*P/\sim$, let $[g,x]$ denote the class of $(g,x)$ in $J$, and let $\pi:G*P\to J$,
$(g,x)\mapsto[g,x]$, be the quotient map. We equip $J$ with the quotient topology. We have $J=\bigsqcup_iJ_i$, where
$J_i:=G\times_iP/\sim$. Indeed, $(g,x)\sim(h,y)$ implies that $gx=hy$; hence, $(g,x),(h,y)\in G\times_iP$ for some $i\in I$. The
map $\pi$ is clearly $I$-graded continuous.

{\unitlength=1bp$$\latexpic(0,0)(-173,36)
\put(0,40){$G*P$}\put(28,44){\vector(1,0){15}}\put(32,47){$\pi$}\put(44,40){$J$}\put(13,35){\vector(1,-2){12}}\put(8,20){$\psi$}
\put(43,35){\vector(-1,-2){12}}\put(41,20){$\varphi$}\put(22,0){$M$}
\endlatexpic$$}

\rightskip=60pt

\vskip-43pt

The groupoid $G$ acts (by homeomorphisms) on $G*P$ by means of the composition~of (composable) arrows: $h(g,x):=(hg,x)$. If
$(g,x)\sim(h,y)$ and $fg$ is defined, then $(fg,x)\sim(fh,y)$. So, we obtain the induced action $h[g,x]:=[hg,x]$ of $G$ on $J$. It
follows from the definition of quotient topology and from $\pi^{-1}(gp)=g\pi^{-1}p$ that $G$ acts on $J$ by homeomorphisms. We get
a commutative diagram of $I$-graded continuous $G$-maps $\psi$, $\pi$, and $\varphi[g,x]:=gx$.

\rightskip0pt

\medskip

Let
$$[P_i]:=\big\{[1_i,x]\mid x\in P_i\big\},\qquad[\overset{\,_\circ}\to P_i]:=\big\{[1_i,x]\mid x\in\overset{\,_\circ}\to
P_i\big\},\qquad[P]:=\bigcup_i[P_i].$$
Clearly, $J=\bigcup\limits_{g,i}g[P_i]$ and $g[\overset{\,_\circ}\to P_i]\cap h[\overset{\,_\circ}\to P_j]\ne\varnothing$ implies
$g=h$. In other words, $[P]$ is a {\it fundamental region\/} for the action of $G$ on $J$.

We assume that $\pi^{-1}[1_i,x]$ is {\bf finite} for every $i\in I$ and $x\in\partial P_i$, hence, for every $x\in P_i$. Take
$x\in P_i$. Let $\pi^{-1}[1_i,x]=\big\{(g_1,x_1),\dots,(g_n,x_n)\big\}$ for some $g_k\in G$ and $x_k\in P_{j_k}$ with
$(g_k,x_k)\in G\times_iP$. The~polyhedra $g_kP_{j_k}$ are the {\it formal neighbours\/} of $P_i$ at $x$. For $\delta>0$, define
$$N_{x_k,\delta}:=\big\{y\in P_{j_k}\mid d(y,x_k)<\delta\big\}\subset P_{j_k},\qquad
N_{x,\delta}:=\bigcup_{k=1}^n(g_k,N_{x_k,\delta})\subset G\times_iP,\qquad W_{x,\delta}:=\pi N_{x,\delta}\subset J_i,$$
where $d(\cdot,\cdot)$ stands for the distance function on $M$. Using this notation, we state the

\medskip

{\bf2.1.~Tessellation condition.} A polyhedron $P$ with a given face-pairing satisfies {\it tessellation condition\/} if

\smallskip

(1) for all $i\in I$ and $x\in P_i$, there exists some $\delta(x)>0$ such that $\pi^{-1}W_{x,\delta}=N_{x,\delta}$ and
$\varphi W_{x,\delta}=B(x,\delta)$ for any $0<\delta\le\delta(x)$;

(2) for every $i\in I$, some open metric neighbourhood $N_i$ of $P_i$ in $M_i$ is tessellated; this means that
$N(P_i,\varepsilon)\subset N_i$ for some $\varepsilon>0$ and that there exists a function $f_i:P_i\to\Bbb R$ taking positive values
such that $\varphi:W_{P_i,f_i}\to N_i$ is bijective, where $W_{P_i,f_i}:=\bigcup\limits_{x\in P_i}W_{x,f_i(x)}\subset J_i$.

\medskip

The proof of the following proposition follows closely that of [AGr2, Proposition 2.2].

\medskip

{\bf2.2.~Proposition.} {\sl Tessellation condition\/ {\rm2.1} implies that\/ $\varphi$ is a regular covering.}

\medskip

{\bf Proof.} First, using only (the finiteness of $\pi^{-1}[1,x]$ and) tessellation condition 2.1 (1), we will show that $\varphi$
is an open map. Indeed, tessellation condition 2.1 (1) and the fact that $G$ acts on $J$ by homeomorphisms imply that
$gW_{x,\delta}$ is open in $J$ for every $(g,x)\in G*P$ and $0<\delta\le\delta(x)$. Moreover,
if~$(g,x)\in G*P$ and $0<\delta\le\delta(x)$, then $\varphi(gW_{x,\delta})=g\varphi W_{x,\delta}=gB(x,\delta)=B(gx,\delta)$ is open
in $M$. Hence, it suffices to prove that $\Cal B:=\big\{gW_{x,\delta}\mid(g,x)\in G*P,\ 0<\delta\le\delta(x)\big\}$ is a base of
the topology on $J$. Let~$[g,x]=g[1_i,x]\in J_j$ and let $U\subset J_j$ be an open neighbourhood of $[g,x]$, where $g:M_i\to M_j$.
The~open set $\pi^{-1}U\subset G\times_jP$ is a union $\pi^{-1}U=\bigcup_h(h,U_h)$, where each $U_h$ is open in some $P_l$,
$h:M_l\to M_j$. Let $\pi^{-1}[1_i,x]:=\big\{(g_1,x_1),\dots,(g_n,x_n)\big\}$. We have
$\pi^{-1}[g,x]=\big\{(gg_1,x_1),\dots,(gg_n,x_n)\big\}$. It follows that $x_k\in U_h$ for $h=gg_k$, $k=1,\dots,n$. Take
$0<\delta\le\delta(x)$ such that $N_{x_k,\delta}\subset U_h$ for $h=gg_k$ and all $k=1,\dots,n$. Clearly,
$gN_{x,\delta}=\bigcup\limits_k(gg_k,N_{x_k,\delta})\subset\pi^{-1}U$. So,
$[g,x]\in gW_{x,\delta}=g\pi N_{x,\delta}=\pi(gN_{x,\delta})\subset\pi\pi^{-1}U=U$ which implies that $\Cal B$ is a base of the
topology on $J$ and that $\varphi$ is open.

\smallskip

Using the fact that $\Cal B$ is a base and that $\pi^{-1}[g,x]$ is finite for every $(g,x)\in G*P$, it is easy to see that $J$ is
Hausdorff.

\smallskip

If $J_i=\varnothing$ for some $i\in I$, then $\varphi|_{J_i}:J_i\to M_i$ is a regular covering of degree zero. So, in order to
prove that $\varphi$ is a regular covering, it suffices to show that $\varphi_i:=\varphi|_{J_i}:J_i\to M_i$ is a surjective regular
covering in the case of $J_i\ne\varnothing$.

\smallskip

We fix such an $i\in I$. Since $\varphi$ is open, $\varphi J_i$ is open in $M_i$. Let $x\in\Cl(\varphi J_i)$. Then
$B(x,\varepsilon)\cap gP_j\ne\varnothing$ for some $g:M_j\to M_i$ and $j\in I$. It follows that
$x\in N(gP_j,\varepsilon)\subset gN_j=\varphi(gW_{P_j,f_j})\subset\varphi J_i$. Hence, $\varphi J_i$ is closed in $M_i$. Since
$M_i$ is connected, $\varphi_i$ is surjective.

Take $x\in M_i$. Define
$$G_x:=\big\{g\in G\mid U_x\cap gP_j\ne\varnothing\;{\text{\rm for }}g:M_j\to M_i,\ j\in I\big\},$$
where $U_x\subset B(x,\frac12\varepsilon)$ is a path-connected open neighbourhood of $x$. For every $G_x\ni g:M_j\to M_i$,~let
$$W_g:=(\varphi_i^{-1}U_x)\cap gW_{P_j,f_j}\subset J_i.$$
Being a continuous open bijection, $\varphi_i:gW_{P_j,f_j}\to gN_j$ is a homeomorphism. Since $U_x\cap gP_j\ne\varnothing$ implies
that $U_x\subset B(x,\frac12\varepsilon)\subset N(gP_j,\varepsilon)\subset gN_j$, we conclude that $\varphi_i:W_g\to U_x$ is a
homeomorphism. Moreover,
$$\varphi_i^{-1}U_x=\bigcup_{g\in G_x}W_g.$$
Indeed, if $\varphi_i[g,y]\in U_x$ with $y\in P_j$ and $g:M_j\to M_i$, then $gy\in U_x$, $g\in G_x$, and $[g,y]\in gW_{P_j,f_j}$,
i.e., $[g,y]\in W_g$.

It remains to show that the distinct $W_g$'s are disjoint. Suppose that $W_{g_1}\cap W_{g_2}\ne\varnothing$ for some
$g_1,g_2\in G_x$. The projection $W_{g_1}\times W_{g_2}\to W_{g_1}$ induces a homeomorphism between
$$X:=\big\{(x_1,x_2)\in W_{g_1}\times W_{g_2}\mid\varphi_ix_1=\varphi_ix_2\big\}$$
and $W_{g_1}$. The diagonal
$$\Delta_{W_{g_1}\cap W_{g_2}}=\Delta_J\cap(W_{g_1}\times W_{g_2})\subset X$$
is closed in $X$ as $J$ is Hausdorff. Therefore, the image $W_{g_1}\cap W_{g_2}$ of $\Delta_{W_{g_1}\cap W_{g_2}}$ is closed in
$W_{g_1}$. Since $W_{g_1}$ is connected, we obtain $W_{g_1}=W_{g_2}$
$_\blacksquare$

\bigskip

\centerline{\bf3.~Compact smooth polyhedra with face-pairing}

\medskip

In this section, we introduce compact smooth polyhedra equipped with face-pairing and observe some local properties of such
polyhedra.

\medskip

{\bf3.1.~Compact smooth polyhedra.} Let $M$ be a (not necessarily connected) smooth manifold. A~compact smooth polyhedron $P$ of
dimension $d=0$ is simply a finite subset $P\subset M$. We define a {\it compact smooth polyhedron\/} $P\subset M$ and its {\it
faces\/} by induction on the dimension $d>0$ of $P$ :

\smallskip

$\bullet$ $P$ is a finite disjoint union $P=\bigsqcup_iF_i$ of its {\it codimension\/} $0$ {\it faces\/} $F_i$. Every $F_i$ is a
connected compact subspace contained in a smooth (connected) submanifold $S_i\subset M$ of dimension $d$, $F_i\subset S_i$. Every
$F_i$ coincides with the closure of its nonempty {\it interior\/} $\varnothing\ne\overset{\,_\circ}\to F_i\subset S_i$, i.e.,
$\overset{\,_\circ}\to F_i$ is open in $S_i$ and $F_i=\Cl\overset{\,_\circ}\to F_i$. We call $S_i$ a {\it submanifold of the
face\/} $F_i$.

$\bullet$ The {\it boundary\/} $\partial F_i:=F_i\setminus\overset{\,_\circ}\to F_i$ of $F_i$ is decomposed into a (possibly empty)
finite union $\partial F_i=\bigcup_js_{ij}$ of connected compact smooth polyhedra of dimension $d-1$, called {\it codimension\/}
$1$ {\it faces\/} of $P$. A codimension $k$ face $f$ of $s_{ij}$ is called a {\it codimension\/} $k+1$ {\it face\/} of $P$. A {\it
submanifold\/} of the face $f$ of $P$ is the same as a submanifold of the face $f$ of $s_{ij}$. For any $j_1\ne j_2$, we require
that $s_{ij_1}\cap s_{ij_2}=\partial s_{ij_1}\cap\partial s_{ij_2}$ and that this intersection is a (possibly empty) union of
entire faces.

\medskip

In a riemannian manifold $M$, denote by $S(x,\delta)$ the sphere of radius $\delta$ centred at $x\in M$.

\medskip

{\bf3.2.~Lemma.} {\sl Let\/ $M$ be a complete riemannian manifold, let\/ $P\subset M$ be a compact smooth polyhedron, and let\/
$x\in P$. Then there exists\/ $\delta(x)>0$ such that\/ $S(x,\delta)\cap P$ is a compact smooth polyhedron for any\/
$0<\delta\le\delta(x)$. We can take\/ $\delta(x)$ such that\/ $S(x,\delta)$ intersects only the faces of\/ $P$ that properly
contain\/ $x$, such that the codimension\/ $k$ faces of\/ $S(x,\delta)\cap P$ are the connected components of\/
$S(x,\delta)\cap f$, where\/ $f$ is a codimension\/ $k$ face\/ of\/ $P$ containing\/ $x$, and such that the interior of a face of\/
$S(x,\delta)\cap P$ comes from the interior of the corresponding face of\/ $P$.}

\medskip

{\bf Proof.} Given a point in a submanifold $x\in S\subset M$, the sphere $S(x,\delta)$ intersects $S$ transversally for all
sufficiently small $\delta$ and, moreover, $S(x,\delta)\cap S$ is a smooth sphere. We apply this fact to submanifolds of all faces
of $P$ and choose $\delta(x)>0$ so small that, for any $0<\delta\le\delta(x)$, the sphere $S(x,\delta)$ intersects only the faces
of $P$ that properly contain $x$. In particular, we can assume that $P$ is connected. By~induction on the dimension of $P$, we
assume that $f_i:=S(x,\delta)\cap s_i$ is a compact smooth polyhedron for any codimension $1$ face $s_i$ of $P$ and
$0<\delta\le\delta(x)$. Denote by $c_{ij}$ the faces of codimension $0$ of $f_i$, $f_i=\bigsqcup_jc_{ij}$. By induction,
$\overset{\,_\circ}\to c_{ij}:=c_{ij}\cap\overset{\,_\circ}\to s_i$. Let $S$ and $S_i$ denote respectively submanifolds of $P$ and
$s_i$. Then $S'_i:=S(x,\delta)\cap S_i$ is a submanifold of $c_{ij}$ and $S':=S(x,\delta)\cap S$ is a sphere containing the compact
set $F:=S(x,\delta)\cap P$.

Take $p\in\overset{\,_\circ}\to c_{ij}$. Note that, in an arbitrarily small open neighbourhood $U$ of $p$, there exists a smooth
simple path $c\subset U\cap F$ beginning at $p$ such that $c\setminus\{p\}\subset S(x,\delta)\cap\overset{\,_\circ}\to P$. Indeed,
we can assume that $U$ does not intersect the other proper faces of $P$ and, by routine arguments of transversality, that $U$ is a
smooth chart of $M$ at $p$ such that the submanifolds $S_i$, $S$, and $S(x,\delta)$ are linear subspaces inside $U$, with $p$ as
the origin. The linear subspace $U\cap S_i$ divides $U\cap S$ into two open half-spaces. One of them $H$ contains a point
$b\in\overset{\,_\circ}\to P$ because $U\cap S_i$ is contained in the closure of $U\cap\overset{\,_\circ}\to P$. It suffices to
show that $S(x,\delta)\cap H\subset\overset{\,_\circ}\to P$. Let $q\in S(x,\delta)\cap H\setminus\overset{\,_\circ}\to P$. As $H$
is convex, the straight segment $[b,q]$ lies in $H$, hence, does not intersect $U\cap S_i$. The first point of $[b,q]$ that does
not belong to $\overset{\,_\circ}\to P$ obviously belongs to the closure of $\overset{\,_\circ}\to P$. Consequently, it belongs to
a proper face of $P$. By the choice of $U$, such a point lies in $s_i$~and, hence, in $U\cap S_i$. A contradiction.

By induction, for any $p\in F$, there exists an arbitrarily short path inside $F$ that begins at $p$ and ends at a point in
$S(x,\delta)\cap\overset{\,_\circ}\to P$.

We will show that every path-connected component $C$ of $F$ is closed in $S(x,\delta)$. Let $q\in\Cl C\setminus C$. Then
$q\in F\setminus\overset{\,_\circ}\to P$, hence, $q$ belongs to some $c_{ij}$. We choose an open path-connected neighbourhood $V$
of $q$ in $S'$ that intersects only those $c_{ij}$'s that contain $q$. Since $q\in\Cl C$, there is a point $p\in V\cap C$. By~the
above, there exist a path inside $V\cap F$ that begins at $p$ and ends at a point in $S(x,\delta)\cap\overset{\,_\circ}\to P$.
So,~we can assume that $p\in V\cap C\cap\overset{\,_\circ}\to P$. There is a simple smooth path $c\subset V$ from $p$ to $q$. The
first point $b$ of $c$ that does not belong to $\overset{\,_\circ}\to P$ belongs to the closure of $\overset{\,_\circ}\to P$.
Hence, $b\in V\cap C\cap c_{ij}$ for some $c_{ij}$. By the choice of $V$, we conclude that $q\in c_{ij}$. Since $c_{ij}$ is
path-connected, we obtain $q\in C$.

We can assume that $C\not\subset\overset{\,_\circ}\to P$ for any path-connected component $C$ of $F$. Otherwise, being open and
closed in the sphere $S'$, the component $C$ has to coincide with $S'$ and we are done. Since every $c_{ij}$ is path-connected, any
path-connected component $C$ of $F$ contains some $c_{ij}$. Therefore, we have finitely many components in $F$.

Moreover, $C=\Cl\overset{\,_\circ}\to C$, where $\overset{\,_\circ}\to C:=C\cap\overset{\,_\circ}\to P$ is open in $S'$. Indeed,
suppose that $C\cap c_{ij}\ne\varnothing$. Then $\overset{\,_\circ}\to c_{ij}\subset c_{ij}\subset C$. By the above, we can find a
simple smooth path $c\subset F$ that begins at an arbitrary point $p\in\overset{\,_\circ}\to c_{ij}$ and such that
$c\setminus\{p\}\subset\overset{\,_\circ}\to C$. So, $C\cap c_{ij}\ne\varnothing$ implies
$\overset{\,_\circ}\to c_{ij}\subset\Cl\overset{\,_\circ}\to C$. Since $\overset{\,_\circ}\to c_{ij}$ is dense in $c_{ij}$, we see
that $C\cap c_{ij}\ne\varnothing$ implies $c_{ij}\subset\Cl\overset{\,_\circ}\to C$
$_\blacksquare$

\medskip

{\bf3.3.~Tesselation at a point.} Let $M$ be a $d$-dimensional smooth manifold.

Let $P_1,P_2\subset M$ be $d$-dimensional compact smooth polyhedra sharing a common codimension $1$ face $s\subset P_1\cap P_2$,
let $p\in s$ be a point in the interior of $s$, and let $p\in B\subset M$ be a neighbourhood of $p$. We~say that $P_1,P_2$ {\it
tessellate\/} $B$ if $B=B_1\cup B_2$ and $B_1\cap B_2=B\cap s$, where $B_i:=B\cap P_i$.

Let $P_i\subset M$, $i=1,\dots,n$, be $d$-dimensional compact smooth polyhedra, let $e\subset\bigcap_iP_i$ be a common codimension
$2$ face of the $P_i$'s, let $p\in e$ be a point in the interior of $e$, and let $p\in B\subset M$ be a neighbourhood of $p$. We
say that the polyhedra $P_i$'s {\it tessellate\/} $B$ if, for every $i$ (the indices are modulo $n$), there is a common codimension
$1$ face $s_i$ of $P_{i-1}$ and $P_i$ containing $e$ such that $B=\bigcup_iB_i$, $B_{i-1}\cap B_i=B\cap s_i$, and
$B_i\cap B_j=B\cap e$ for all $i,j$ with $|i-j|>1$, where $B_i:=B\cap P_i$.

\medskip

{\bf3.4.~Remark.} {\sl Let\/ $M$ be a\/ $d$-dimensional complete riemannian manifold, let\/ $P_i\subset M$ be\/ $d$-dimensio\-nal
compact smooth polyhedra sharing a common codimension\/ $1$ or\/ $2$ face\/ $f\subset\bigcap_iP_i$, and let\/ $p\in f$ be a point
in the interior of\/ $f$. Suppose that the\/ $P_i$'s tessellate a neighbourhood of\/ $p$. Let\/ $S(x,\delta)$ be a sphere related
to every polyhedra\/ $P_i$ as in Lemma\/ {\rm3.2} and let\/ $p\in S(x,\delta)$. Then the polyhedra\/ $S(x,\delta)\cap P_i$
tessellate a neighbourhood of\/ $p$ in\/ $S(x,\delta)$}
$_\blacksquare$

\medskip

{\bf3.5.~Smooth polyhedra with face-pairing.} Let $M$ be an oriented riemannian manifold with $d$-dimensional components and let
$P\subset M$ be a $d$-dimensional compact smooth polyhedron. Suppose that we are given a face-pairing, i.e., an involution
$\overline{\phantom{m}}:S\to S$ on the set $S$ of all codimension $1$ faces of $P$ and isometries $I_s\in\Isom M$ satisfying
$I_ss=\overline s$ and $I_{\overline s}=I_s^{-1}$ for all $s\in S$. If every codimension $2$ face $e$ of $P$ belongs to exactly two
codimension $1$ faces $s_1,s_2$ of $P$ (in symbols: $s_1\diamond e\diamond s_2$) and each $I_s$ maps any face of $s$ onto a face of
$\overline s$, we say that $P$ is equipped with {\it face-pairing.}

\medskip

{\bf3.6.~Geometric cycles.} Let $P\subset M$ be a compact smooth polyhedron equipped with face-pairing. Start with
$\overline s_0\diamond e\diamond s_1$. Applying $I_{s_1}$ to $s_1$ and $e$, we obtain $\overline s_1\diamond I_{s_1}e\diamond s_2$.
Applying $I_{s_2}$ to $s_2$ and $I_{s_1}e$, we obtain $\overline s_2\diamond I_{s_2}I_{s_1}e\diamond s_3$, and so on. Since the
number of faces is finite, we eventually arrive back at $\overline s_0\diamond e\diamond s_1$.

A cyclic sequence
$$\overset{^{I_{s_n}}}\to\longrightarrow\overline s_n=\overline s_0\diamond e\diamond
s_1\overset{^{I_{s_1}}}\to\longrightarrow\overline s_1\diamond I_{s_1}e\diamond s_2\overset{^{I_{s_2}}}\to\longrightarrow\overline
s_2\diamond I_{s_2}I_{s_1}e\diamond s_3\overset{^{I_{s_3}}}\to\longrightarrow\dots\overset{^{I_{s_{n-1}}}}
\to\longrightarrow\overline s_{n-1}\diamond I_{s_{n-1}}\cdots I_{s_1}e\diamond s_n\overset{^{I_{s_n}}}\to\longrightarrow,$$
where each term is obtained from the previous one by the above rule, is called a {\it cycle of codimension\/ $2$ faces.} The number
$n$ is the {\it length\/} of the cycle and the isometry $I:=I_{s_n}\cdots I_{s_1}$ will be referred to as the {\it cycle isometry.}
A cycle can be read backwards, i.e., in {\it opposite orientation,} which inverts its isometry. If the cycle isometry is the
identity and if the cycle is the shortest one with this property, then the cycle is said to be {\it geometric.} Clearly, every
cycle is a multiple of a shortest, combinatorial one. Note that, in~a geometric cycle, a codimension $2$ face $e$ may occur several
times in the form $s\diamond e\diamond s'$, where $s,s'$ are codimension $1$ faces containing $e$ (this does not happen to a
combinatorial cycle).

\medskip

{\bf3.7.~Formal neighbours in codimension $\le2$.} Let $P$ be a compact smooth polyhedron with face-pairing.

It is easy to describe all formal neighbours (see section 2 for the definition) of $P$ at a point $x\in s$ in the interior of a
codimension $1$ face $s$. Since $s$ is a unique codimension $1$ face containing $x$ and $\overline s$  is a unique codimension $1$
face containing $I_sx$, we have $\pi^{-1}[1,x]=\big\{(1,x),(I_{\overline s},I_sx)\big\}$. Hence, the only formal neighbours of $P$
at $x$ are $P,I_{\overline s}P$.

Let $x\in e$ be a point in the interior of a codimension $2$ face $e$ and suppose that $e$ belongs to a geometric cycle
$$\overset{^{I_{s_n}}}\to\longrightarrow\overline s_n=\overline s_0\diamond I_0e\diamond
s_1\overset{^{I_{s_1}}}\to\longrightarrow\overline s_1\diamond I_1e\diamond s_2\overset{^{I_{s_2}}}\to\longrightarrow\overline
s_2\diamond I_2e\diamond s_3\overset{^{I_{s_3}}}\to\longrightarrow\dots\overset{^{I_{s_{n-1}}}}\to\longrightarrow\overline
s_{n-1}\diamond I_{n-1}e\diamond s_n\overset{^{I_{s_n}}}\to\longrightarrow,\leqno{\bold{(3.8)}}$$
where $I_j:=I_{s_j}\cdots I_{s_1}$ for $j=0,1,\dots,n$ (we consider $j$ modulo $n$). Then
$\pi^{-1}[1,x]=\big\{(I_j^{-1},I_jx)\mid j=0,1,\dots,n-1\big\}$ and the $I_j^{-1}P$ are all the formal neighbours of $P$ at $x$.
Indeed, suppose that $(I_j^{-1},I_jx)\sim(h,y)$. We can assume that $I_jx\in s'$, $I_{s'}I_jx=y$, and $h^{-1}I_j^{-1}=I_{s'}$ for
some $s'\in S$. In particular, the codimension $1$ face $s'$ intersects the interior of the codimension $2$ face $I_je$. By 3.1, we
have $I_je\diamond s'$. By 3.5, either $s'=\overline s_j$ or $s'=s_{j+1}$. Therefore, either $(h,y)=(I_{j-1}^{-1},I_{j-1}x)$ or
$(h,y)=(I_{j+1}^{-1},I_{j+1}x)$. It remains to observe that the $I_j$, $j=0,1,\dots,n-1$, are all distinct because we could
otherwise take a shorter cycle whose isometry would be the identity.

\medskip

{\bf3.9.~Interior into exterior condition.} Let $P$ be a compact smooth polyhedron equipped with face-pairing. Assuming that $M$ is
oriented, we orient every codimension $1$ face $s$ in such a way that the unit normal vector $n_s$ to $s$ points toward the
interior of $P$. We say that the face-pairing isometries of $P$ send {\it interior into exterior\/} if $I_sn_s=-n_{\overline s}$
for every codimension $1$ face $s\in S$. Since every codimension $1$ face is connected, it suffices to verify this condition at a
single point of every codimension $1$ face.

\medskip

{\bf3.10.~Lemma.} {\sl Let\/ $P$ be a compact smooth polyhedron with face-pairing. Then the face-pairing isometries of\/ $P$ send
interior into exterior iff, for every interior point\/ $x\in\overset{\,_\circ}\to s$ in any codimension\/ $1$ face\/~$s$, the
polyhedra\/ $P,I_{\overline s}P$ tessellate a neighbourhood of\/ $x$.}

\medskip

{\bf Proof.} The proof follows [AGr2, First step of the proof of Theorem 3.5]. Let $x\in s$ be an interior point in a codimension
$1$ face of $P$. We choose $\delta>0$ with the following properties: $B(x,\delta)$ does not intersect the faces of $s$,
$B(I_sx,\delta)$ does not intersect the faces of $\overline s$, $B(x,\delta)\cap\partial P=B(x,\delta)\cap s$,
and~$B(I_sx,\delta)\cap\partial P=B(I_sx,\delta)\cap\overline s$. Let $B:=B(x,\delta)$, $B_1:=B\cap P$, and
$B_2:=B\cap I_{\overline s}P=I_{\overline s}\big(B(I_sx,\delta)\cap P\big)$. Note that $B\cap s\subset B_1\cap B_2$.

Assume that the face-pairing isometries send interior into exterior. Then $B_1\ne B_2$. Pick a point $q_0\in B_1\setminus B_2$ such
that $q_0\notin s$. By the choice of $\delta$, a smooth oriented curve $\gamma\subset B$ connecting $q_0$ and $q\in B\setminus s$
can intersect $\partial P$ and $\partial I_{\overline s}P$ only along $(s\setminus\partial s)\cap B$. We can take such
intersections as being transversal. Due to the interior into exterior condition, when intersecting $s$, the curve $\gamma$ leaves
$B_1$ and enters $B_2$ or {\it vice-versa.} Hence, $q$ belongs to exactly one of $B_1$ and $B_2$. This implies that $B=B_1\cup B_2$
and $B_1\cap B_2=B\cap s$.

The converse is immediate
$_\blacksquare$

\medskip

{\bf3.11.~Total angle condition.} Let $P$ be a compact smooth polyhedron whose face-pairing isometries send interior into exterior.
Pick a point $x\in e$ in a codimension $2$ face $e$ and denote by $\N_xe:=(\T_xe)^\perp$ the plane normal to $e$ at $x$. By 3.5, we
have $\overline s_0\diamond e\diamond s_1$ for suitable codimension $1$ faces $\overline s_0,s_1$. Denote by $n_0,n_1$ the unit
normal vectors to $\overline s_0,s_1$ at $x$ that point towards the interior of $P$. Let $t_0\in\T_x\overline s_0\cap\N_xe$
and\break

\vskip-3pt

\noindent
\hskip341pt$\vcenter{\hbox{\epsfbox{Picture.eps}}}$

\rightskip118pt

\vskip-56pt

\noindent
$t_1\in\T_xs_1\cap\N_xe$ stand for the unit vectors that point respectively towards the interiors of $\overline s_0$ and $s_1$. The
basis $t_0,n_0$ orients $\N_xe$. The {\it interior angle\/} $\alpha_0$ between $\overline s_0$ and $s_1$ at $x$ is the angle from
$t_0$ to $t_1$. Such an angle takes values in $[0,2\pi]$. Note that interchanging the faces $\overline s_0,s_1$ in the definition
does not alter $\alpha_0$.

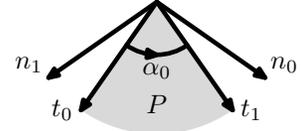
\rightskip0pt

Suppose that $e$ belongs to a geometric cycle (3.8). Denote by $\alpha_j$ the interior angle between $\overline s_j$ and $s_{j+1}$
at $I_jx$. The sum $\sum_{j=0}^{n-1}\alpha_j$ is the {\it total interior angle\/} of the cycle at $x$. It is easy to see that
altering the orientation of the cycle keeps the same values of the $\alpha_j$'s. Using the facts that the face-pairing isometries
send interior into exterior and that the cycle (3.8) is geometric, as in [AGr2, Subsection 3.3], we conclude that
$\alpha:=\sum_{j=0}^{n-1}\alpha_j\equiv0\mod2\pi$. The {\it total angle condition\/} consists in the requirement that
$\alpha=2\pi$. Since any codimension $2$ face is connected, the total angle condition for a single point $x$ implies the same
condition for all interior points of the codimension $2$ faces involved in the geometric cycle.

\medskip

{\bf3.12.~Lemma.} {\sl Let\/ $P$ be a compact smooth polyhedron whose face-pairing isometries send interior into exterior, let\/
$e$ be a codimension\/ $2$ face participating in a geometric cycle\/ {\rm(3.8),} and let\/ $x\in\overset{\,_\circ}\to e$ be an
interior point. Then the total angle of the cycle at\/ $x$ equals\/ $2\pi$ iff the polyhedra\/ $I_j^{-1}P$ tessellate a
neighbourhood of\/ $x$.}

\medskip

{\bf Proof.} The proof is similar to that of Lemma 3.10 (see also [AGr2, First step of the proof of Theorem 3.5]). We choose
$\delta>0$ with the following properties: $B(I_jx,\delta)$ does not intersect the faces of $I_je$, $B(I_jx,\delta)$ does not
intersect any face of $\overline s_j$ or $s_{j+1}$ except $I_je$, and
$B(I_jx,\delta)\cap\partial P=\big(B(I_jx,\delta)\cap\overline s_j\big)\cup\big(B(I_jx,\delta)\cap s_{j+1}\big)$ for all
$j=0,1,\dots,n-1$. Let $B:=B(x,\delta)$ and $B_j:=B\cap I_j^{-1}P=I_j^{-1}\big(B(I_jx,\delta)\cap P\big)$. Note that
$B\cap I_j^{-1}s_{j+1}\subset B_j\cap B_{j+1}$ and that $B\cap e\subset B_i\cap B_j$ for all $i,j$. Let $T_j\subset\N_xe$ stand for
the closed sector containing the interior angle of $I_j^{-1}P$ at $x$ and let $F:=B\cap\bigcup_jI_j^{-1}s_{j+1}$.

Assume that the total angle of the cycle at $x$ equals $2\pi$. Since the face-pairing isometries send interior into exterior, we
have $\bigcup_jT_j=N_xe$ and $\overset{\circ}\to T_i\cap\overset{\circ}\to T_j=\varnothing$ if $i\not\equiv j\mod n$. Hence,
$B_i\ne B_j$ for $i\not\equiv j\mod n$. Pick a point $q_0$ that lies in exactly one of the $B_j\setminus F$. By the choice of
$\delta$, a smooth oriented curve $\gamma\subset B$ connecting $q_0$ and $q\in B\setminus F$ may intersect
$\bigcup_j\partial I_j^{-1}P$ only along $F$. We can assume that $\gamma$ does not intersect $e$ and is transversal to $F$. Due to
the interior into exterior condition, when intersecting $I_j^{-1}s_{j+1}$, the curve $\gamma$ leaves $B_j$ and enters $B_{j+1}$ or
{\it vice-versa.} Hence, $q$ belongs to exactly one of the $B_j$'s. It follows that $B=\bigcup_jB_j$ and, by the description of
formal neighbours in 3.7, that $B_j\cap B_{j+1}=B\cap I_j^{-1}s_{j+1}$ and $B_i\cap B_j=B\cap e$ for all $i,j$ with $|i-j|>1$.

Conversely, assume that the polyhedra $I_j^{-1}P$ tessellate a neighbourhood of $x$. We can assume (dimi\-nishing $\delta$, if
needed) that they tessellate $B$. Then $\bigcup_jT_j=N_xe$ because $B=\bigcup_jB_j$. Finally, we have
$\overset{\circ}\to T_i\cap\overset{\circ}\to T_j=\varnothing$ if $i\not\equiv j\mod n$ since
$B_j\cap B_{j+1}=B\cap I_j^{-1}s_{j+1}$ and $B_i\cap B_j=B\cap e$ for all $i,j$ with $|i-j|>1$. In other words, the total angle of
the cycle at $x$ equals $2\pi$
$_\blacksquare$

\medskip

As a corollary to the proofs of Lemmas 3.10 and 3.12, the tessellation by formal neighbours at points in the interior of
codimension $1$ or $2$ faces implies, for such points, a condition that is stronger than tessellation condition 2.1 (1).

\medskip

{\bf3.13.~Corollary.} {\sl Let\/ $P$ be a compact smooth polyhedron such that every codimension\/ $2$ face participates in a
geometric cycle. Assume that, for every point\/ $x$ in the interior of a codimension\/ $1$ or\/ $2$ face, the formal neighbours
of\/ $P$ at\/ $x$ tessellate a neighbourhood of\/ $x$. Then, for every such\/ $x$, there exists $\delta(x)>0$ such that
tessellation condition\/ {\rm2.1 (1)} holds and\/ $\varphi:W_{x,\delta}\to B(x,\delta)$ is a bijection for all\/
$0<\delta\le\delta(x)$.}

\medskip

{\bf Proof.} Let $x$ be in the interior of a codimension $1$ face $s$. We take $\delta(x)>0$ such that the open ball
$B\big(x,\delta(x)\big)$ is tessellated by the formal neighbours $P,I_{\overline s}P$ of $P$ at $x$ (see 3.7). Let
$0<\delta\le\delta(x)$. As in the proof of Lemma 3.10, we can assume that $B(x,\delta)$ does not intersect the faces of $s$,
$B(I_sx,\delta)$ does not intersect the faces of $\overline s$, $B(x,\delta)\cap\partial P=B(x,\delta)\cap s$, and
$B(I_sx,\delta)\cap\partial P=B(I_sx,\delta)\cap\overline s$. Let~$B:=B(x,\delta)$, $B_1:=B\cap P$, and
$B_2:=B\cap I_{\overline s}P$. It follows from $B=B_1\cup B_2$ and $B_1\cap B_2=B\cap s$ that $\varphi:W_{x,\delta}\to B(x,\delta)$
is a bijection. The choice of $\delta(x)$ and the description of formal neighbours given in 3.7 imply that
$\pi^{-1}W_{x,\delta}=N_{x,\delta}$.

Let $x$ be in the interior of a codimension $2$ face $e$ participating in a geometric cycle (3.8). We choose $\delta(x)>0$ such
that the open ball $B\big(x,\delta(x)\big)$ is tessellated by the formal neighbours $I_j^{-1}P$ of $P$ at $x$ (see~3.7). Let
$0<\delta\le\delta(x)$. We assume that $B(x,\delta)$ satisfies the properties in the beginning of the proof of Lemma 3.12:
$B(I_jx,\delta)$ does not intersect the faces of $I_je$, $B(I_jx,\delta)$ does not intersect any face of $\overline s_j$ or
$s_{j+1}$ except $I_je$, and
$B(I_jx,\delta)\cap\partial P=\big(B(I_jx,\delta)\cap\overline s_j\big)\cup\big(B(I_jx,\delta)\cap s_{j+1}\big)$ for all
$j=0,1,\dots,n-1$. Let $B:=B(x,\delta)$ and $B_j:=B\cap I_j^{-1}P$. Since $B=\bigcup_jB_j$,
$B_j\cap B_{j+1}=B\cap I_j^{-1}s_{j+1}$, and $B_i\cap B_j=B\cap e$ for all $i,j$ with $|i-j|>1$, we conclude that
$\varphi:W_{x,\delta}\to B(x,\delta)$ is a bijection. The choice of $\delta(x)$ and the description of formal neighbours given in
3.7 imply the rest
$_\blacksquare$

\medskip

The above strengthened tessellation condition 2.1 (1) implies
tessellation condition 2.1 (2) :

\medskip

{\bf3.14.~Lemma.} {\sl Let\/ $P$ be a compact smooth polyhedron with
face-pairing in a complete riemannian manifold. Assume that\/ $P$
satisfies tessellation condition\/ {\rm2.1 (1)} and that\/
$\varphi:W_{x,\delta}\to B(x,\delta)$ is a bijection for every\/
$x\in P$ and\/ $0<\delta\le\delta(x)$. Then\/ $P$ satisfies
tessellation condition\/ {\rm 2.1.}}

\medskip

{\bf Proof.} As in section 2, we put
$[P_i]:=\big\{[1_i,x]\mid x\in P_i\big\}$. Note that
$\varphi:W_{x,\delta}\to B(x,\delta)$ is a homeomorphism for every
$x\in P_i$ and $0<\delta\le\delta(x)$ since $\varphi$ is an open map
(see the very beginning of the proof of Proposition 2.2). In view of
the commutative diagram in section 2, the map $\varphi:[P_i]\to P_i$ is
bijective. Due to the compactness of $P_i$, there exists a family of
open neighbourhoods $N_n$ of $[P_i]$ whose closures are compact,
$\Cl N_{n+1}\subset N_n$ for all $n$, and $\bigcap_nN_n=[P_i]$.

Suppose that $\varphi|_{N_n}$ is not injective for every $n\in\Bbb N$.
Hence, there are sequences $\{p_n\},\{q_n\}$ such that
$p_n,q_n\in N_n$, $p_n\ne q_n$, and $\varphi p_n=\varphi q_n$ for every
$n\in\Bbb N$. We can assume that the sequences converge. Clearly, their
limits belong to $[P_i]$. As the map $\varphi:[P_i]\to P_i$ is
bijective, the limits coincide. This contradicts the fact that
$\varphi$ is injective in a neighbourhood $W_{p,\delta}$ of the limit
$[1_i,p]$. Therefore, $\varphi:N_n\to\varphi N_n$ is a homeomorphism
for some $n\in\Bbb N$ which implies tessellation condition 2.1 (2) by
the compactness of $P_i$
$_\blacksquare$

\medskip

{\bf3.15.~Definition.} A compact smooth polyhedron $P\subset M$ is said to be {\it connected near codimension\/}~$3$~if, for every
$v\in F$, there exists $\delta(v)>0$ such that $B(v,\delta)\cap P\setminus F$ is connected for all $0<\delta<\delta(v)$, where $F$
stands for the union of all codimension $3$ faces of $P$.

\medskip

The following example illustrates that the interior into exterior and total angle conditions do not imply tessellation condition
2.1 in the absence of connectedness near codimension $3$.

\medskip

{\bf3.16.~Example.} Let $\Bbb S^2\subset\Bbb H_\Bbb R^3$ be a metric sphere centred at a point $s\in\Bbb H_\Bbb R^3$ in real
hyperbolic space and let $a,b,c\in\Bbb S^2$ be vertices of a right-angled triangle $T_1\subset\Bbb S^2$. Denote by $m$ the middle
point of the side $bc$. Let $T_2$ stand for the triangle $T_1$ reflected in $s$. The edges and vertices of the polygon
$T:=T_1\cup T_2\subset\Bbb S^2$ are the edges and vertices of the triangles and the faces of $T$ are the two disjoint triangles.
Face-pairing is given by the rotation about $a$ by $\pi/2$ and the rotation about $m$ by $\pi$. The polyhedron
$P\subset\Bbb H_\Bbb R^3$ is the cone of those rays emanated from $s$ to the absolute that pass through points of~$T$. Extend
correspondingly the isometries. One readily verifies the exterior into exterior and total angle conditions and sees that the copies
of $P$ with respect to the (finite) group generated by the face-pairing isometries tessellate {\it twice\/} the space
$\Bbb H_\Bbb R^3$
$_\blacksquare$

\medskip

The idea of the example is that the formally tessellated space $J$ (see section 2) is a pair of disjoint balls with their centres
identified. Even in dimension $3$, it seems possible to construct a more sophisticated example where the tessellation around
vertices provides a similar picture but with different numbers of glued balls at different vertices.

\medskip

{\bf3.17.~Lemma.} {\sl Let\/ $P\subset M$ be a connected near codimension\/ $3$ compact smooth polyhedron in a riemannian
manifold\/ $M$ of dimension\/ $\le4$ and let\/ $v\in F$ be a point in the union\/ $F$ of all codimension\/ $3$ faces of\/ $P$. Then
there exists\/ $\delta_0>0$ such that\/ $S(v,\delta)\cap P\setminus F$ is connected for all\/ $0<\delta<\delta_0$. In~particular,
$S(v,\delta)\cap P$ is connected for such\/ $\delta$.}

\medskip

{\bf Proof.} For sufficiently small $\delta$, the sphere $S(v,\delta)$ intersects only the faces of $P$ that properly contain $v$
and is transversal to the submanifolds of these faces. It follows that $B(v,\delta_0)\cap P\setminus F\to(0,\delta_0)$ is a
(locally) trivial bundle with the fibres $S(v,\delta)\cap P\setminus F$. The rest follows from the fact that the finite set
$S(v,\delta)\cap F$ lies in the closure of $S(v,\delta)\cap P\setminus F$
$_\blacksquare$

\medskip

The concept of connectedness near codimension $3$ has a purely topological nature. For instance, a~compact smooth polyhedron
$P\subset M$ in a $3$-manifold $M$ is not connected near codimension $3$ iff there exist a codimension $3$ face $v\in P$ and an
open neighbourhood $v\in U\subset M$ such that the set $U\cap P\setminus\{v\}$ possesses path-connected components $C_1,C_2$ with
$v\in\Cl C_i$ for $i=1,2$. Indeed, suppose that, for some codimension $3$ face $v\in P$, the set $B(v,\delta)\cap P\setminus\{v\}$
is not connected for an arbitrarily small $\delta>0$. As in the proof of Lemma 3.17, we have a trivial bundle
$B(v,\delta_0)\cap P\to(0,\delta_0)$ for some $\delta_0>0$. So,~the fibre $S(v,\delta)\cap P$ possesses nonempty path-connected
components $A_1,A_2\subset S(v,\delta)\cap P$. It remains to put $U:=B(v,\delta_0)$ and $C_i:=A_i\times(0,\delta_0)$ in the sense
of the identification $B(v,\delta_0)\cap P\setminus\{v\}\simeq\big(S(v,\delta)\cap P\big)\times(0,\delta_0)$. The converse is
obvious.

\medskip

{\bf3.18.~Lemma.} {\sl Let\/ $P\subset M$ be a connected near codimension\/ $3$ compact smooth polyhedron in a complete
riemannian\/ $4$-manifold\/ $M$ and let\/ $v\in F$ be a point in the union\/ $F$ of all codimension\/ $3$ faces of\/ $P$. Then
there exists\/ $\delta_0>0$ such that the compact smooth polyhedron\footnote{We need the completeness of $M$ in order to be able to
apply Lemma 3.2 thus guaranteeing that $S(v,\delta)\cap P$ is a compact smooth polyhedron for small $\delta$.}\/
$S(v,\delta)\cap P$ in the\/ $3$-manifold\/ $S(v,\delta)$ is connected near codimension $3$ for any $0<\delta<\delta_0$.}

\medskip

{\bf Proof.} As in the proof of Lemma 3.17, we have a trivial bundle $B(v,\delta_0)\cap P\to(0,\delta_0)$ with the fibres
$S(v,\delta)\cap P$. Suppose that $S(v,\delta)\cap P$ is not connected near codimension $3$ for some $0<\delta<\delta_0$. Then
there exist a codimension $3$ face $f$ of $P$, a point $p\in S(v,\delta)\cap f$, an open neighbourhood $p\in U\subset S(v,\delta)$
of $p$ in $S(v,\delta)$, and path-connected components $C_1,C_2$ of $U\cap P\setminus\{p\}$ such that $p\in\Cl C_i$ for $i=1,2$.
In~terms of a suitable identification $B(v,\delta_0)\cap P\setminus\{v\}\simeq\big(S(v,\delta)\cap P\big)\times(0,\delta_0)$ such
that $\{p\}\times(0,\delta_0)\subset f$, we obtain the path-connected components $\widehat C_i:=C_i\times(0,\delta_0)$ of
$\big(U\times(0,\delta_0)\big)\cap P\setminus\big(\{p\}\times(0,\delta_0)\big)$ such that
$f\supset\{p\}\times(0,\delta_0)\subset\Cl\widehat C_i$ for $i=1,2$. Taking an open ball
$B(p,\varepsilon)\subset U\times(0,\delta_0)$, $\varepsilon>0$, we arrive at a contradiction with connectedness of $P$ near
codimension $3$
$_\blacksquare$

\bigskip

\centerline{\bf4.~Main theorems}

\medskip

{\bf4.1.~Theorem.} {\sl Let\/ $M$ be an oriented complete riemannian manifold with\/ $3$-dimensional components and let\/
$P\subset M$ be a compact smooth polyhedron connected near codimension $3$. Suppose that the face-pairing isometries of\/ $P$ send
interior into exterior, that every codimension\/ $2$ face participates in a geometric cycle, and that, for every such cycle, the
total interior angle equals\/ $2\pi$ at some point of a codimension\/ $2$ face of the cycle. Then tessellation condition\/ {\rm2.1}
is satisfied.}

\medskip

{\bf Proof.} Let $v_1,\dots,v_m$ denote the  codimension $3$ faces of $P$. By combining Lemmas 3.10 and 3.12 with Corollary 3.13
and Lemma 3.14, we reduce the theorem to the following facts:

\smallskip

$\bullet$ For every $i$, there are finitely many formal neighbours of $P$ at $v_i$.

$\bullet$ There exists $\delta_0>0$ such that $\pi^{-1}W_{v_i,\delta}=N_{v_i,\delta}$ and
$\varphi:W_{v_i,\delta}\to B(v_i,\delta)$ is a bijection for all $i$ and $0<\delta\le\delta_0$.

\medskip

{\bf4.2.1.~Polyhedron $P_\delta\subset M_\delta$ and its face-pairing.} By Lemma 3.2, there exists a small $\delta_0>0$ such that
$P_\delta:=\bigsqcup_iS(v_i,\delta)\cap P$ is a compact smooth polyhedron in the smooth $2$-manifold
$M_\delta:=\bigsqcup_iS(v_i,\delta)$ for any $0<\delta\le\delta_0$. Moreover, for any $i$, the $2$-sphere $S(v_i,\delta)$
intersects only those faces of $P$ that properly contain $v_i$. The codimension $k$ faces of $P_\delta$ are the connected
components $(v_i,f)_j$ of $S(v_i,\delta)\cap f$, where $f\subset P$ is a codimension $k$ face of $P$ containing $v_i$. The interior
of $(v_i,f)_j$ is the intersection of $(v_i,f)_j$ with the interior of $f$. The only possible faces of $P_\delta$ are those of
codimensions $0,1,2$. By~Lemma~3.17, we assume that $S(v_i,\delta)\cap P$ is connected for all $i$ and $0<\delta\le\delta_0$. In
other words, $S(v_i,\delta)\cap P$ are the codimension $0$ faces of $P_\delta$.

We equip $P_\delta$ with face-pairing as follows. Let $(v_i,s)_j$ be a codimension $1$ face of $P_\delta$, where $s\subset P$ is a
codimension $1$ face of $P$ and $v_i\in s$. Since $I_s$ maps $v_i$ to $I_sv_i$, $S(v_i,\delta)$ onto $S(I_sv_i,\delta)$, and $s$
onto $\overline s$, it maps the component $(v_i,s)_j$ of $S(v_i,\delta)\cap s$ onto a certain component $(I_sv_i,\overline s)_l$ of
$S(I_sv_i,\delta)\cap\overline s$. We~define the involution $\overline{\phantom{m}}:S_\delta\to S_\delta$ on the set $S_\delta$ of
codimension $1$ faces of $P_\delta$ as $\overline{(v_i,s)}_j:=(I_sv_i,\overline s)_l$. The isometry
$I_{(v_i,s)_j}\in\Isom M_\delta$ is the restriction of $I_s\in\Isom M$.

Let $G_\delta\subset\Isom M_\delta$ denote the subcategory generated by the face-pairing isometries of $P_\delta$.

\vskip10pt

{\unitlength=1bp$$\latexpic(0,0)(224,50)
\put(0,54){$G_\delta*P_\delta$}\put(33,57){\vector(1,0){31}}\put(44,60){$\beta$}\put(66,54){$G*P$}\put(17,51){\vector(0,-1){42}}
\put(6,29){$\psi_\delta$}\put(23,51){\vector(2,-3){11}}\put(29,44){$\pi_\delta$}\put(72,51){\vector(-2,-3){10}}\put(61,44){$\pi$}
\put(79,51){\vector(0,-1){42}}\put(82,29){$\psi$}\put(33,27){$J_\delta$}\put(44,30){\vector(1,0){10}}\put(44,33){$\gamma$}
\put(55,27){$J$}\put(33,25){\vector(-2,-3){11}}\put(31,15){$\varphi_\delta$}\put(61,26){\vector(2,-3){12}}\put(58,15){$\varphi$}
\put(11,0){$M_\delta$}\put(25,4){\vector(1,0){47}}\put(73,0){$M$}
\endlatexpic$$}

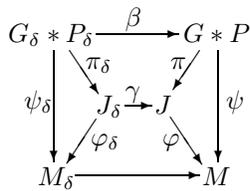
\leftskip102pt

\vskip-48pt

{\bf4.2.2.~Diagram.} The polyhedron $P$ and the subcategory $G\subset\Isom M$ generated by face-pairing isometries determine the
commutative triangle (see section~2) on the right side of the diagram. Analogously, we obtain from $P_\delta$ and $G_\delta$ the
commutative triangle on the left side of the diagram.

The face-pairing isometries that generate $G_\delta$ are restrictions of the face-pairing isometries that generate $G$. Sending the
component $S(v_i,\delta)$ of $M_\delta$ to the component of $M$ that contains $S(v_i,\delta)$ and using the fact that the
composition of restric-\break

\vskip-12pt

\leftskip0pt

\noindent
tions is the restriction of the composition, we can define a functor $\alpha:G_\delta\to G$. At the level of generators,
$\alpha:I_{(v_i,s)_j}\mapsto I_s$. In order to verify the correctness of the definition of $\alpha$ for morphisms, we need only to
show that an isometry $I:M_j\to M_j$ of a component $M_j$ of $M$ containing $S(v_i,\delta)$ is the identity if its restriction
$I|_{S(v_i,\delta)}$ is the identity. This can be achieved by a suitable choice of $\delta_0$. For sufficiently small $\delta_0$,
the sphere $S(v_i,\delta)$ has a unique centre with respect to the radius $\delta$. Hence, $Iv_i=v_i$. If $\delta_0$ is less than
the radius of injectivity of $M$ at $v_i$, then $I=1$. Indeed, let $p\in M_j$ and let $\Gamma$ be a shortest geodesic joining $v_i$
and $p$. If necessary, we can extend the geodesic $\Gamma$ so that it intersects $S(v_i,\delta)$, say, at $x\in S(v_i,\delta)$.
Since $Iv_i=v_i$, $Ix=x$, and $v_i,x\in\Gamma$, we obtain $I|_\Gamma=1$ and, therefore, $Ip=p$.

We define a map $\beta:G_\delta*P_\delta\to G*P$, $\beta:(g,x)\mapsto(\alpha g,x)$, in the commutative exterior square of the
diagram (we think of $P_\delta$ as being included into $P$). Note that $\beta:G_\delta*P_\delta\to G*P$ is injective. Indeed,
$\beta(g,x)=\beta(h,y)$ implies $x=y$, $\alpha g=\alpha h$, and $gx=(\alpha g)x=(\alpha h)y=hy$. So, if $x\in S(v_i,\delta)$ and
$gx\in S(v_j,\delta)$, then $g,h:S(v_i,\delta)\to S(v_j,\delta)$, $g=\alpha g|_{S(v_i,\delta)}$, and $h=\alpha h|_{S(v_i,\delta)}$,
implying the fact.

The equivalence relation $\sim_\delta$ on $G_\delta*P_\delta$ is induced by the equivalence relation $\sim$ on $G*P$, i.e.,
$(g,x)\sim_\delta(h,y)$ iff $\beta(g,x)\sim\beta(h,y)$ for any $(g,x),(h,y)\in G_\delta*P_\delta$. Indeed, for $x\in(v_i,s)_j$, the
relation $(I_{(v_i,s)_j},x)\sim_\delta(1,I_{(v_i,s)_j}x)$ is nothing but the relation $(I_s,x)\sim(1,I_sx)$. The rest follows from
$\beta\big(g(h,x)\big)=(\alpha g)\beta(h,x)$. We obtain an injective map $\gamma:J_\delta\to J$ in the commutative top square of
the diagram of continuous maps. Since $\pi_\delta$ is surjective, the entire diagram is commutative.

\medskip

{\bf4.2.3.~Formal neighbours in $M_\delta$ and in $M$.} Let $x$ be a point in a codimension $1$ face $(v_i,s)_j$ of
$S(v_i,\delta)\cap P_\delta$ and let $\pi^{-1}[1,x]=\big\{(g_1,x_1),\dots,(g_n,x_n)\big\}$. Using the description of formal
neighbours given in 3.7, we are going to show that there exist $h_k\in G_\delta$ such that $\beta(h_k,x_k)=(g_k,x_k)$ for all
$k=1,\dots,n$.

Suppose that $x$ belongs to the interior of $s$. Then $\pi^{-1}[1,x]=\big\{(1,x),(I_s^{-1},I_sx)\big\}$. The point $I_sx$ belongs
to a component of $S(I_sv_i,\delta)\cap\overline s$. Taking $h_k:=I_{(v_i,s)_j}^{-1}$, we have $\beta(h_k,I_sx)=(I_s^{-1},I_sx)$.

Suppose that $x$ lies in the interior of a codimension $2$ face $e\subset s$ of $P$. The face $e$ belongs to a geometric cycle and,
using the notation from (3.8) with $s_1:=s$, we have $\pi^{-1}[1,x]=\big\{(I_k^{-1},I_kx)\mid k=0,1,\dots,n-1\big\}$. By induction,
we get $h_{k-1}\in G_\delta$ such that $\beta(h_{k-1},I_{k-1}x)=(I_{k-1}^{-1},I_{k-1}x)$ and the point
$I_{k-1}x\in S(I_{k-1}v_i,\delta)\cap I_{k-1}e\subset S(I_{k-1}v_i,\delta)\cap s_k$ belongs to a component of
$S(I_{k-1}v_i,\delta)\cap s_k$. The isometry $I_{s_k}$ maps this component onto a certain component of
$S(I_kv_i,\delta)\cap\overline s_k$. Therefore, the~corresponding restriction $h^{-1}\in G_\delta$ of $I_{s_k}$ provides
$\beta(h_{k-1}h,I_kx)=(I_{k-1}^{-1}I_{s_k}^{-1},I_{s_k}I_{k-1}x)=(I_k^{-1},I_kx)$.

In both cases, $(h_k,x_k)\in\pi_\delta^{-1}[1_i,x]$ because $\gamma$ is injective. Conversely, let
$(h,y)\in\pi_\delta^{-1}[1_i,x]$. Then $\pi\beta(h,y)=\gamma[1_i,x]=[1,x]$ and $(h,y)=(h_k,x_k)$ for some $k$ as $\beta$ is
injective. In other words, there is a one-to-one correspondence between formal neighbours of $P$ at $x$ and formal neighbours of
$P_\delta$ at $x$.

At the level of $M$, this correspondence means that the formal neighbours of $P_\delta$ at $x\in(v_i,s)_j$ are the intersections
with $S(v_i,\delta)$ of the corresponding formal neighbours of $P$ at $x$. Indeed, given $(g_k,x_k)\in\pi^{-1}[1,x]$ and
$(h_k,x_k)\in G_\delta*P_\delta$ as above, let $S(v_{i_k},\delta)\subset M_{j_k}$ denote the components of $M_\delta$ and $M$ that
contain $x_k$. Since $g_kx_k=x$, we obtain $g_kS(v_{i_k},\delta)=S(v_i,\delta)$. By the choice of $\delta_0$, we have
$S(v_{i_k},\delta)\cap P_\delta=S(v_{i_k},\delta)\cap P_{j_k}$. Hence,
$\varphi_\delta\big[h_k,S(v_{i_k},\delta)\cap P_\delta\big]=\varphi\big[g_k,S(v_{i_k},\delta)\cap P_{j_k}\big]$ by the
commutativity of the bottom square. We obtain
$h_k\big(S(v_{i_k},\delta)\cap P_\delta\big)=g_k\big(S(v_{i_k},\delta)\cap P_{j_k}\big)=S(v_i,\delta)\cap g_kP_{j_k}$.

\medskip

{\bf4.2.4.~Tessellation of $M_\delta$.} Let $x\in P_\delta$. The above bijection between $\pi_\delta^{-1}[1,x]$ and
$\pi^{-1}[1,x]$ implies that $\pi_\delta^{-1}[1,x]$ is finite. Combining Lemma 3.10, Lemma 3.12, and Remark 3.4 with the fact
obtained in 4.2.3 that the formal neighbours of $P_\delta$ at any $x\in P_\delta$ are the intersections of the formal neighbours of
$P$ with $S(v_i,\delta)\ni x$, we conclude that formal neighbours of $P_\delta$ at any $x\in P_\delta$ tessellate a neighbourhood
of $x$.

Every codimension $2$ face $x\in S(v_i,\delta)$ of $P_\delta$ participates in a geometric cycle of $P_\delta$. Indeed, $x$ lies in
the interior of a codimension $2$ face $e$ of $P$ and $e$ belongs to a geometric cycle (3.8) of $P$. As in 4.2.3, we~can `restrict'
this cycle to a cycle of $P_\delta$. The isometry of the cycle of $P_\delta$ is the identity. The~isometry of a shorter cycle of
$P_\delta$ cannot be the identity because the corresponding polyhedra $I_{\delta,j}^{-1}P_\delta$ tessellate a neighbourhood of $x$
in $M_\delta$ and, hence, the interiors of these polyhedra are disjoint in this neighbourhood, which makes impossible the equality
$I_{\delta,j}=1$ for $0<j<n$.

By Corollary 3.13, Lemma 3.14, and Proposition 2.2, $\varphi_\delta:J_\delta\to M_\delta$ is a surjective regular covering.

\medskip

{\bf4.2.5.}~Let $C_\delta$ be a component of $J_\delta$. If $(g_j,p)\in\pi_\delta^{-1}C_\delta$ and $p\in F_j$ belongs to a
codimension $0$ face $F_j$ of $P_\delta$, then $(g_j,F_j)\subset\pi_\delta^{-1}C_\delta$ because $F_j$ is path-connected. Moreover,
$\pi_\delta^{-1}C_\delta$ is a finite union of such pieces, $\pi_\delta^{-1}C_\delta=\bigsqcup\limits_{j=1}^k(g_j,F_j)$. Indeed,
the measure of (the interior of) a codimension $0$ face of $P_\delta$ in $M_\delta$ is limited from below by a positive constant,
$\varphi_\delta:C_\delta\to S(v_i,\delta)$ is a homeomorphism for some $i$, and
$\pi_\delta(g_j,\overset{\,_\circ}\to F_j)\subset C_\delta$ are all disjoint for $(g_j,F_j)\subset\pi_\delta^{-1}C_\delta$.

Note that $\pi_\delta^{-1}C_\delta$ is a minimal nonempty subset of the type $\bigsqcup\limits_{j=1}^k(g_j,F_j)$ closed with
respect to taking $\sim_\delta$-equivalent elements. Otherwise, $C_\delta$ is a disjoint union of two nonempty compacts.

\medskip

{\bf4.2.6.}~To any $(g,v)\in\pi^{-1}[1,v_i]$ and $0<\delta\le\delta_0$, we associate
$$N_{(g,v),\delta}:=B(v,\delta)\cap P,\qquad F_{(g,v),\delta}:=S(v,\delta)\cap P\subset P_\delta.$$
Since the codimension $0$ face of $P$ containing $v$ is path-connected, by Lemma 3.17, $F_{(g,v),\delta}$ is a (nonempty)
codimension $0$ face of $P_\delta$ for a suitable choice of $\delta_0$. Let
$$N_{v_i,\delta}:=\bigsqcup_{(g,v)\in\pi^{-1}[1,v_i]}(g,N_{(g,v),\delta})\subset G*P,\qquad
S_{v_i,\delta}:=\bigsqcup_{(g,v)\in\pi^{-1}[1,v_i]}(g,F_{(g,v),\delta})\subset\beta(G_\delta*P_\delta).$$
The inclusions easily follow from $(g,v)\in G*P$ and from a decomposition of $g$ in generators $I_s$'s due to the choice of
$\delta_0$. Obviously, $\psi S_{v_i,\delta}\subset S(v_i,\delta)$. We have $(g,F_{(g,v),\delta})=\beta(g',F_{(g,v),\delta})$ for
some $g'\in G_\delta$.

Since $N_{v_i,\delta}=\pi^{-1}[1,v_i]\sqcup\bigsqcup\limits_{0<\delta'<\delta}S_{v_i,\delta'}$, in order to show that
$\pi^{-1}(\pi N_{v_i,\delta})=N_{v_i,\delta}$, it suffices to observe that $\pi^{-1}(\pi S_{v_i,\delta})=S_{v_i,\delta}$, i.e.,
that $S_{v_i,\delta}$ is closed with respect to taking $\sim$-equivalent elements. Let $(g,I_sx)\in S_{v_i,\delta}$, where
$x\in s$. We need to show that $(gI_s,x)\in S_{v_i,\delta}$. For some $v_j$, we have $(g,v_j)\in\pi^{-1}[1,v_i]$ and
$I_sx\in S(v_j,\delta)\cap\overline s$. Since $v_j\in\overline s$ by the choice of $\delta_0$, we conclude that
$(gI_s,I_{\overline s}v_j)\sim(g,v_j)\in\pi^{-1}[1,v_i]$ and $x\in S(I_{\overline s}v_j,\delta)\cap P$. So,
$x\in F_{(gI_s,I_{\overline s}v_j),\delta}$ and $(gI_s,x)\in S_{v_i,\delta}$.

Moreover, $\pi_\delta S_\delta$ is path-connected, where
$S_\delta:=\bigsqcup\limits_{(g,v)\in\pi^{-1}[1,v_i]}(g',F_{(g,v),\delta})$. Indeed, $F_{(g,v),\delta}$ is path-connected for any
$(g,v)\in\pi^{-1}[1,v_i]$ by Lemma 3.17 and by the choice of $\delta_0$. Let $v_j\in s$. Consider $F_{(g,v_j),\delta}$,
$F_{(gI_{\overline s},I_sv_j),\delta}$, and the corresponding $g',g''\in G_\delta$. Then $\pi_\delta(g',F_{(g,v_j),\delta})$ and
$\pi_\delta(g'',F_{(gI_{\overline s},I_sv_j),\delta})$ have a common point because $S(v_j,\delta)\cap s\ne\varnothing$ due to the
path-connectedness of $s$ and to the choice of~$\delta_0$. In particular, $\pi_\delta S_\delta$ lies in a component $C_\delta$ of
$J_\delta$ homeomorphic to $S(v_i,\delta)$, $\varphi_\delta:C_\delta\to S(v_i,\delta)$.

On the other hand, $S_\delta$ is closed with respect to taking $\sim_\delta$-equivalent elements because
$S_{v_i,\delta}=\beta S_\delta$ is closed with respect to taking $\sim$-equivalent elements and $\sim_\delta$ is induced by $\sim$.
By 4.2.5, $\pi_\delta S_\delta=C_\delta$ for a suitable component $C_\delta$ of $J_\delta$ such that
$\varphi_\delta:C_\delta\to S(v_i,\delta)$ is a bijection.

Summarizing, $\varphi:\pi N_{v_i,\delta}\to B(v_i,\delta)$ is a bijection and $\pi^{-1}(\pi N_{v_i,\delta})=N_{v_i,\delta}$. Since
$F_{(g,v),\delta}\ne\varnothing$ for every $(g,v)\in\pi^{-1}[1,v_i]$, the finiteness of $\pi^{-1}[1,v_i]$ follows from 4.2.5
$_\blacksquare$

\medskip

{\bf4.3.~Theorem.} {\sl Let\/ $M$ be an oriented complete riemannian manifold with\/ $4$-dimensional components and let\/ $P$ be a
compact smooth polyhedron connected near codimension $3$. Suppose that the face-pairing isometries of\/ $P$ send interior into
exterior, that every codimension\/ $2$ face participates in a geometric cycle, and that, for every such cycle, the total interior
angle equals\/ $2\pi$ at some point of a codimension\/ $2$ face of the cycle. If the stabilizer of every codimension\/ $3$ face in
the groupoid generated by face-pairing isometries induces on the face a finite group of isometries,\footnote{This condition is
obviously satisfied for codimension $3$ faces that are not round circles.}
then tessellation condition\/ {\rm2.1} is satisfied.}

\medskip

{\bf Proof} mostly follows the line of the proof of Theorem 4.1.

\medskip

{\bf4.4.1.~Extra codimension $4$ faces.} First, we introduce additional codimension $4$ faces inside codimension $3$ faces in order
to make the latter topological segments (not circles).

Let $f_i$ be a codimension $3$ face of $P$ inside a codimension $1$ face $s_i$ of $P$, $f_i\subset s_i$ for $i=1,2$. We say that
$f_1,f_2$ are equivalent if $s_2=\overline s_1$ and $f_2=I_{s_1}f_1$ and take the closure of this symmetric relation with respect
to transitivity.

Let $f$ be a face in such an equivalence class. By conditions of Theorem 4.3, the group $\Gamma$ of isometries of $f$ induced by
the stabilizer of $f$ in the groupoid is finite. Pick an interior point $p\in\overset{\,_\circ}\to f$ and introduce the new
codimension $4$ faces $\Gamma p\subset f$. Then, using face-pairing isometries, we copy these new faces into the other codimension
$3$ faces equivalent to $f$.

\medskip

{\bf4.4.2.}~By combining Lemmas 3.10 and 3.12 with Corollary 3.13 and Lemma 3.14, we reduce the theorem to the following facts:

\smallskip

$\bullet$ For every $p$ in a codimension $3$ face of $P$, there are finitely many formal neighbours of $P$ at $p$.

$\bullet$ For every $p$ in a codimension $3$ face of $P$, there exists $\delta(p)>0$ such that $\pi^{-1}W_{p,\delta}=N_{p,\delta}$
and $\varphi:W_{p,\delta}\to B(p,\delta)$ is a bijection for all $0<\delta\le\delta(p)$.

\medskip

{\bf4.4.3.~Polyhedron $P_{\varepsilon,\delta}\subset M_{\varepsilon,\delta}$ and its face-pairing.} Fix some
$\varepsilon\in[0,\frac12]$. In every codimension $3$ face $f\subset P$, we take the points that are on the distance
$\varepsilon\ell$, measured along $f$, from one of the vertices of $f$, where $\ell$ stands for the length of $f$. So, we take a
couple of points in every codimension $3$ face $f$ if $\varepsilon\ne\frac12$ and the middle point of $f$, otherwise. Denote by
$v_1,\dots,v_m$ all such points ($\varepsilon$ is fixed). By~Lemma~3.2, there exists a small $\delta_0>0$ such that
$P_{\varepsilon,\delta}:=\bigsqcup_iS(v_i,\delta)\cap P$ is a compact smooth polyhedron in the smooth $3$-manifold
$M_{\varepsilon,\delta}:=\bigsqcup_iS(v_i,\delta)$ for any $0<\delta\le\delta_0$. Moreover, for any $i$, the~$3$-sphere
$S(v_i,\delta)$ intersects only those faces of $P$ that properly contain $v_i$. The codimension $k$ faces of $P_\delta$ are the
connected components $(v_i,f)_j$ of $S(v_i,\delta)\cap f$, where $f\subset P$ is a codimension $k$ face of $P$ containing $v_i$.
The interior of $(v_i,f)_j$ is the intersection of $(v_i,f)_j$ with the interior of $f$. The only possible faces of
$P_{\varepsilon,\delta}$ are those of codimensions $0,1,2,3$. By Lemma 3.17, we assume that $S(v_i,\delta)\cap P$ are connected for
all $i$ and $0<\delta\le\delta_0$. In other words, $S(v_i,\delta)\cap P$ is the codimension $0$ faces of $P_{\varepsilon,\delta}$.

We equip $P_{\varepsilon,\delta}$ with face-pairing exactly as in 4.2.1 and denote by
$G_{\varepsilon,\delta}\subset\Isom M_{\varepsilon,\delta}$ the subcategory generated by face-pairing isometries of
$P_{\varepsilon,\delta}$.

\medskip

{\bf4.4.4.~Diagram.} After adapting notation, this subsection is literally the same as subsection 4.2.2.

\medskip

{\bf4.4.5.~Formal neighbours in $M_{\varepsilon,\delta}$ and in $M$ in codimension $\le2$.} Let $x\in S(v_i,\delta)\cap P$ be not a
codimension $3$ face of $P_{\varepsilon,\delta}$. Then, literally following 4.2.3, we show that the formal neighbours of
$P_{\varepsilon,\delta}$ at $x$ are the intersections with $S(v_i,\delta)$ of the corresponding formal neighbours of $P$ at $x$.

\medskip

{\bf4.4.6.~Tesselation of $M_{\varepsilon,\delta}$.} By Lemma 3.18, the polyhedron
$P_{\varepsilon,\delta}\subset M_{\varepsilon,\delta}$ is connected near codimension $3$. By 4.4.5, Lemma 3.10, Lemma 3.12, and
Remark 3.4, the polyhedron $P_{\varepsilon,\delta}$ satisfies the conditions of Theorem 4.1 because the geometric cycles of $P$
being `restricted' to $P_{\varepsilon,\delta}$ remain geometric. Indeed, the isometry of a shorter cycle cannot be the identity for
the same reason as in 4.2.4. So,~by~Theorem~4.1 and Proposition 2.2,
$\varphi_{\varepsilon,\delta}:J_{\varepsilon,\delta}\to M_{\varepsilon,\delta}$ is a surjective regular covering.

\medskip

{\bf4.4.7.}~As in 4.2.5, we show that $\pi_{\varepsilon,\delta}^{-1}C_{\varepsilon,\delta}=\bigsqcup\limits_{j=1}^k(g_j,F_j)$ is a
minimal nonempty subset of this type closed with respect to taking $\sim_{\varepsilon,\delta}$-equivalent elements for any
component $C_{\varepsilon,\delta}$ of $J_{\varepsilon,\delta}$, where the $F_j$'s are codimension $0$ faces of
$P_{\varepsilon,\delta}$.

\medskip

{\bf4.4.8.}~Here, after adapting notation, we follow literally 4.2.6 using 4.4.7 in place of 4.2.5
$_\blacksquare$

\bigskip

\centerline{\bf5.~References}

\medskip

[AGr1] S.~Anan$'$in, C.~H.~Grossi, {\it Coordinate-free classic geometries,} Mosc.~Math.~J.~{\bf11} (2011), 633--655, see also
http://arxiv.org/abs/math/0702714

[AGr2] S.~Anan$'$in, C.~H.~Grossi, {\it Yet another Poincar\'e polyhedron theorem,} Proc.~Edinb.~Math.~Soc. {\bf54} (2011),
297--308, see also http://arxiv.org/abs/0812.4161

[Ana] S.~Anan$'$in, {\it Research plan for\/} 2012--2013

[CaS] D.~I.~Cartwright, T.~Steger, {\it Enumeration of the\/ $50$ fake projective planes,} C.~R.~Math.~{\bf348} (2010), 11--13

[Mas] B.~Maskit, {\it Kleinian Groups,} Grundlehren der mathematichen Wissenschaften {\bf287}, Springer-Verlag, 1988. xiii+326 pp.

[Mos] G.~D.~Mostow, {\it On a remarkable class of polyhedra in complex hyperbolic space,} Pacific J.~Math. {\bf86} (1980), 171--276

[Mum] D.~Mumford, {\it An algebraic surface with\/ $K$ ample, $(K^2)=9$, $p_g=q=0$,} Amer.~J.~Math.~{\bf101} (1979), 233--244

[Rei] I.~Reider, {\it Geography and the number of moduli of surfaces of general type,} Asian J.~Math.~{\bf9} (2005), no.~3,
407--448

[Yau] S.~T.~Yau, {\it Calabi's conjecture and some new results in algebraic geometry,} Proc.~Natl.~Acad.~Sci. {\bf74} (1977),
1798--1799

\enddocument